\documentclass[twosided,reqno]{amsart}

\usepackage{amsmath}
\usepackage{amssymb}
\usepackage{amsthm}

\theoremstyle{theorem}
\newtheorem{theorem}{\scshape Theorem }[section]

\theoremstyle{definition}

\numberwithin{equation}{section}

\begin{document}

\title[higher-order Bernoulli, Euler and Hermite polynomials]{Some identities of higher-order Bernoulli, Euler and Hermite polynomials arising from umbral calculus}

\author{Dae San Kim$^1$}
\address{$^1$ Department of Mathematics, Sogang University, Seoul 121-742, Republic of Korea.}
\email{dskim@sogang.ac.kr}
\author{Taekyun Kim$^2$}
\address{$^2$ Department of Mathematics, Kwangwoon University, Seoul 139-701, Republic of Korea.}
\email{tkkim@kw.ac.kr}
\author{Dmitry V. Dolgy $^3$}
\address{$^3$ Hanrimwon, Kwangwoon University, Seoul 139-701, Republic of Korea.}
\email{dgekw2011@gmail.com}

\author{Seog-Hoon Rim$^4$}
\address{$^4$ Department of Mathematics Education, Kyungpook National University, Taegu 702-701, Republic of Korea.}
\email{shrim@knu.ac.kr}


\maketitle

\begin{abstract}
In this paper, we study umbral calculus to have alternative ways of obtaining our results. That is, we derive some interesting identities of the higher-order Bernoulli, Euler and Hermite polynomials arising from umbral calculus to have alternative ways.
\end{abstract}

\section{Introduction}
As is well known, the Hermite polynomials are defined by the generating function to be
\begin{equation}\label{1}
e^{2xt-t^2}=e^{H(x)t}=\sum_{n=0} ^{\infty} H_n (x) \frac{t^n}{n!},
\end{equation}
with the usual convention about replacing $H^n(x)$ by $H_n(x)$ (see \cite{08, 10}).
In the special case, $x=0$, $H_n(0)=H_n$ are called the {\it{$n$-th Hermite numbers}}. The  Bernoulli polynomials of order $r$ are given by the generating function to be
\begin{equation}\label{2}
\left(\frac{t}{e^t-1}\right)^re^{xt}=\sum_{n=0} ^{\infty} B_n ^{(r)} (x) \frac{t^n}{n!},~(r \in {\mathbb{R}}).
\end{equation}
From \eqref{2}, the $n$-th Bernoulli numbers of order $r$ are defined by $B_n ^{(r)}(0)=B_n^{(r)}$ (see [1-16]). The higher-order Euler polynomials are also defined by the generating function to be
\begin{equation}\label{3}
\left(\frac{2}{e^t+1}\right)^re^{xt}=\sum_{n=0} ^{\infty}E_n ^{(r)} (x)\frac{t^n}{n!},~(r \in {\mathbb{R}}),
\end{equation}
and $E_n ^{(r)}(0)=E_n ^{(r)}$ are called the {\it{$n$-th Euler numbers}} of order $r$ (see [1-16]).

The first Stirling number is given by
\begin{equation}\label{4}
(x)_n=x(x-1)\cdots(x-n+1)=\sum_{l=0} ^n S_1(n,k)x^l,{\text{ (see [6,13])}},
\end{equation}
and the second Stirling number is defined by the generating function to be
\begin{equation}\label{5}
(e^t-1)^n=n!\sum_{l=n} ^{\infty} S_2(l,n) \frac{t^l}{l!},{\text{ (see [6,9,13])}}.
\end{equation}
For $\lambda(\neq 1) \in {\mathbb{C}}$, the {\it{Frobenius-Euler polynomials}} are given by
\begin{equation}\label{6}
\left(\frac{1-\lambda}{e^t-\lambda}\right)^re^{xt}=\sum_{n=0} ^{\infty} H_n ^{(r)}(x|\lambda)\frac{t^n}{n!},~(r \in {\mathbb{R}}){\text{ (see [1,5])}}.
\end{equation}
In the special case, $x=0$, $H_n ^{(r)}(0|\lambda)=H_n ^{(r)}(\lambda)$ are called the {\it{$n$-th Frobenius-Euler numbers}} of order $r$.

Let ${\mathcal{F}}$ be the set of all formal power series in the variable $t$ over ${\mathbb{C}}$ with
\begin{equation*}
{\mathcal{F}}=\left\{ \left.f(t)=\sum_{k=0} ^{\infty} \frac{a_k}{k!} t^k~\right|~ a_k \in {\mathbb{C}} \right\}.
\end{equation*}
Let us assume that ${\mathbb{P}}$ is the algebra of polynomials in the variable $x$ over ${\mathbb{C}}$ and ${\mathbb{P}}^{*}$ is the vector space of all linear functionals on ${\mathbb{P}}$. $\left< L~|~p(x)\right>$ denotes the action of the linear functional $L$ on a polynomial $p(x)$ and we remind that the vector space structure on ${\mathbb{P}}^{*}$ is defined by
\begin{equation*}
\begin{split}
\left<L+M|p(x) \right>&=\left<L|p(x) \right>+\left<M|p(x) \right> ,\\
\left<cL|p(x) \right>&=c\left<L|p(x) \right>,
\end{split}
\end{equation*}
where $c$ is a complex constant (see \cite{06,09,13}).

The formal power series $f(t)=\sum_{k=0} ^{\infty} \frac{a_k}{k!}t^k \in {\mathcal{F}}$ defines a linear functional on ${\mathbb{P}}$ by setting
\begin{equation}\label{7}
\left<f(t)|x^n \right>=a_n,{\text{ for all }} n\geq 0{\text{ (see [6,9,13])}}.
\end{equation}
Then, by \eqref{7}, we get
\begin{equation}\label{8}
\left<t^k | x^n \right>=n! \delta_{n,k},~(n,k \geq 0),
\end{equation}
where $\delta_{n,k}$ is the Kronecker symbol (see \cite{06,09,13}).

Let $f_L(t)=\sum_{k=0} ^{\infty} \frac{\left<L|x^k\right>}{k!}t^k$ (see \cite{13}). For $f_L(t)=\sum_{k=0} ^{\infty} \frac{\left<L|x^k\right>}{k!}t^k$, we have $\left<f_L(t)|x^n\right>=\left<L|x^n\right>$. The map $L \mapsto f_L(t)$ is a vector space isomorphism from ${\mathbb{P}}^{*}$ onto ${\mathcal{F}}$. Henceforth, ${\mathcal{F}}$ will be thought of as both a formal power series and a linear functional. We shall call ${\mathcal{F}}$ the {\it{umbral algebra}}. The umbral calculus is the study of umbral algebra (see \cite{06,09,13}).

The order $o(f(t))$ of the non-zero power series $f(t)$ is the smallest integer $k$ for which the coefficient of $t^k$ does not vanish. A series $f(t)$ having $o(f(t))=1$ is called a {\it{delta series}} and a series $f(t)$ having $o(f(t))=0$ is called an {\it{invertible series}}. Let $f(t)$ be a delta series and $g(t)$ be an invertible series. Then there exists a unique sequence $S_n(x)$ of polynomials such that $\left<g(t)f(t)^k|S_n(x)\right>=n!\delta_{n,k}$ where $n,k \geq 0$. The sequence $S_n(x)$ is called {\it{Sheffer sequence}} for $(g(t),f(t))$, which is denoted by $S_n(x)\sim (g(t),f(t))$. By \eqref{7} and \eqref{8}, we see that $\left<e^{yt}|p(x)\right>=p(y)$. For $f(t)\in {\mathcal{F}}$ and $p(x) \in {\mathbb{P}}$, we have
\begin{equation}\label{9}
f(t)=\sum_{k=0} ^{\infty} \frac{\left<f(t)|x^k\right>}{k!} t^k,~p(x)=\sum_{k=0} ^{\infty} \frac{\left<t^k|p(x)\right>}{k!}x^k,
\end{equation}
and, by \eqref{9}, get
\begin{equation}\label{10}
p^{(k)}(0)=\left<t^k|p(x)\right>,~\left<1|p^{(k)}(x)\right>=p^{(k)}(0).
\end{equation}
Thus, from \eqref{10}, we have
\begin{equation}\label{11}
t^kp(x)=p^{(k)}(x)=\frac{d^kp(x)}{dx^k}.
\end{equation}
In \cite{06,09,13}, we note that $\left<f(t)g(t)|p(x)\right>=\left<g(t)|f(t)p(x)\right>$.

For $S_n(x) \sim \left(g(t),f(t)\right)$, we have
\begin{equation}\label{11}
\frac{1}{g({\bar{f}}(t))}e^{y{\bar{f}}(t)}=\sum_{k=0} ^{\infty} \frac{S_k(y)}{k!}t^k,{\text{ for all }}y \in {\mathbb{C}},
\end{equation}
where ${\bar{f}}(t)$ is the compositional inverse of $f(t)$. For $S_n(x) \sim \left(g(t),f(t)\right)$ and $r_n(x)=(h(t),l(t))$, let us assume that
\begin{equation}\label{12}
S_n(x)=\sum_{k=0} ^n C_{n,k} r_k(x),{\text{ (see [6,9,13])}}.
\end{equation}
Then we have
\begin{equation}\label{13}
C_{n,k}=\frac{1}{k!}\left. \left< \frac{h({\bar{f}}(t))}{g({\bar{f}}(t))}l({\bar{f}}(t))^k\right| x^n \right>,{\text{ (see [13])}}.
\end{equation}
The equation \eqref{12} and \eqref{13} are called the alternative ways of Sheffer sequences.

In this paper, we study umbral calculus to have alternative ways of obtaining our results. That is, we derive some interesting identities of the higher-order Bernoulli, Euler and Hermite polynomials arising from umbral calculus to have alternative ways.

\section{Some identities of higher-order Bernoulli, Euler and Hermite polynomials}

In this section, we use umbral calculus to have alternative ways of obtaining our results. Let us consider the following Sheffer sequences:
\begin{equation}\label{14}
E_n ^{(r)}(x)\sim\left(\left(\frac{e^t+1}{2}\right)^r,t\right),H_n(x)\sim\left(e^{\frac{1}{4}t^2},\frac{t}{2}\right).
\end{equation}
Then, by \eqref{12}, we assume that
\begin{equation}\label{15}
E_n ^{(r)}(x)=\sum_{k=0} ^n C_{n,k}H_k(x).
\end{equation}
From \eqref{13} and \eqref{15}, we have
\begin{equation}\label{16}
\begin{split}
C_{n,k}&=\frac{1}{k!}\left. \left<\frac{e^{\frac{1}{4}t^2}}{\left(\frac{e^t+1}{2}\right)^r}\left(\frac{t}{2}\right)^k \right| x^n\right> \\
&=\frac{1}{k!2^k} \left. \left< \left(\frac{2}{e^t+1}\right)^re^{\frac{1}{4}t^2}\right|t^kx^n\right> \\
&=2^{-k}\binom{n}{k}\left.\left<\left(\frac{2}{e^t+1}\right)^r\right| e^{\frac{1}{4}t^2}x^{n-k}\right> \\
&=2^{-k}\binom{n}{k}\left< \left.\left(\frac{2}{e^t+1}\right)^r\right| \sum_{l=0} ^{\left[\frac{n-k}{2}\right]} \frac{1}{4^ll!}t^{2l}x^{n-k}\right> \\
&=2^{-k}\binom{n}{k}\sum_{l=0} ^{\left[\frac{n-k}{2}\right]} \frac{1}{2^{2l}l!}(n-k)_{2l}\left<\left.1\right. \left|\left(\frac{2}{e^t+1}\right)^rx^{n-k-2l}\right. \right>\\
&=2^{-k}\binom{n}{k}\sum_{l=0} ^{\left[\frac{n-k}{2}\right]} \frac{1}{2^{2l}l!}(n-k)_{2l} E_{n-k-2l} ^{(r)} \\
&=n!\sum_{0 \leq l \leq n-k,~l:{\text{even}}}\frac{E_{n-k-l} ^{(r)}}{\left(\frac{l}{2}\right)!2^{k+l}k!(n-k-l)!}.
\end{split}
\end{equation}
Therefore, by \eqref{15} and \eqref{16}, we obtain the following theorem.
\begin{theorem}
For $n \geq 0$, we have
\begin{equation*}
E_n ^{(r)} (x)=n!\sum_{k=0} ^n \left\{\sum_{0 \leq l \leq n-k,~l:{\text{even}}}\frac{E_{n-k-l} ^{(r)}}{k!(n-k-l)!2^{k+l}\left(\frac{l}{2}\right)!}\right\}H_k(x).
\end{equation*}
\end{theorem}
Let us consider the following two Sheffer sequences:
\begin{equation}\label{17}
B_n ^{(r)}(x)\sim\left(\left(\frac{e^t-1}{t}\right)^r,t\right),~H_n(x)\sim\left(e^{\frac{1}{4}t^2},\frac{t}{2}\right).
\end{equation}
Let us assume that
\begin{equation}\label{18}
B_n ^{(r)}(x)=\sum_{k=0} ^n C_{n,k}H_k(x).
\end{equation}
By \eqref{13} and \eqref{17}, we get
\begin{equation}\label{19}
\begin{split}
C_{n,k}&=\frac{1}{k!}\left. \left<\frac{e^{\frac{1}{4}t^2}}{\left(\frac{e^t-1}{t}\right)^r}\left(\frac{t}{2}\right)^k \right| x^n\right> \\
&=2^{-k}\binom{n}{k}\left.\left<\left(\frac{t}{e^t-1}\right)^r\right| \sum_{l=0} ^{\infty} \left(\frac{1}{4}\right)^l \frac{1}{l!} t^{2l}x^{n-k}\right> \\
&=2^{-k}\binom{n}{k}\sum_{l=0} ^{\left[\frac{n-k}{2}\right]} \frac{1}{l!4^l}(n-k)_{2l} \left. \left<\left(\frac{t}{e^t-1}\right)^r \right|x^{n-k-2l}\right>\\
&=2^{-k}\binom{n}{k}\sum_{l=0} ^{\left[\frac{n-k}{2}\right]} \frac{(n-k)!}{l!2^{2l}(n-k-2l)!}\left< \left.1 \right. \left| \left(\frac{t}{e^t-1} \right)^rx^{n-k-2l}\right. \right> \\
&=2^{-k}\binom{n}{k}\sum_{l=0} ^{\left[\frac{n-k}{2}\right]} \frac{(n-k)!}{l!2^{2l}(n-k-2l)!}B_{n-k-2l} ^{(r)} \\
&=n!\sum_{0 \leq l \leq n-k,~l:{\text{even}}}\frac{B_{n-k-l}  ^{(r)}}{k!(n-k-l)!2^{k+l}\left(\frac{l}{2}\right)!}.
\end{split}
\end{equation}
Therefore, by \eqref{18} and \eqref{19}, we obtain the following theorem.
\begin{theorem}
For $n \geq 0$, we have
\begin{equation*}
B_n ^{(r)}(x)=n!\sum_{k=0} ^n \left\{\sum_{0 \leq l \leq n-k,~l:{\text{even}}}\frac{B_{n-k-l}  ^{(r)}}{k!(n-k-l)!2^{k+l}\left(\frac{l}{2}\right)!} \right\}H_k(x).
\end{equation*}
\end{theorem}
Consider
\begin{equation}\label{20}
H_n ^{(r)}(x|\lambda)\sim\left(\left(\frac{e^t-\lambda}{1-\lambda}\right)^r,t\right),H_n(x)\sim\left(e^{\frac{1}{4}t^2},\frac{t}{2}\right).
\end{equation}
Let us assume that
\begin{equation}\label{21}
H_n ^{(r)}(x|\lambda)=\sum_{k=0} ^n C_{n,k}H_k(x).
\end{equation}
By \eqref{13}, we get
\begin{equation}\label{22}
\begin{split}
C_{n,k}&=\frac{1}{k!}\left. \left<\frac{e^{\frac{1}{4}t^2}}{\left(\frac{e^t-\lambda}{1-\lambda}\right)^r}\left(\frac{t}{2}\right)^k \right| x^n\right> \\
&=\frac{1}{k!2^k}\left.\left<\left(\frac{1-\lambda}{e^t-\lambda}\right)^re^{\frac{1}{4}t^2}\right| t^k x^n\right> \\
&=2^{-k}\binom{n}{k}\sum_{l=0} ^{\left[\frac{n-k}{2}\right]} \frac{(n-k)_{2l}}{l!4^l}\left< \left.1 \right. \left| \left(\frac{1-\lambda}{e^t-\lambda} \right)^rx^{n-k-2l}\right. \right> \\
&=n!\sum_{l=0} ^{\left[\frac{n-k}{2}\right]}\frac{H_{n-k-2l}  ^{(r)}(\lambda)}{l!2^{2l+k}(n-k-2l)!k!} \\
&=n!\sum_{0 \leq l \leq n-k,~l:{\text{even}}}\frac{H_{n-k-l}  ^{(r)}(\lambda)}{\left(\frac{l}{2}\right)!2^{k+l}(n-k-l)!k!}.
\end{split}
\end{equation}
Therefore, by \eqref{21} and \eqref{22}, we obtain the following theorem.
\begin{theorem}
For $n \geq 0$, we have
\begin{equation*}
H_n ^{(r)} (x|\lambda)=n!\sum_{k=0} ^n \left\{\sum_{0 \leq l \leq n-k,~l:{\text{even}}}\frac{H_{n-k-l}  ^{(r)}(\lambda)}{k!(n-k-l)!\left(\frac{l}{2}\right)!2^{k+l}}\right\}H_k(x).
\end{equation*}
\end{theorem}
Let us assume that
\begin{equation}\label{23}
H_n (x)=\sum_{k=0} ^n C_{n,k}E_k ^{(r)}(x).
\end{equation}
From \eqref{13}, \eqref{14} and \eqref{23}, we have
\begin{equation}\label{24}
\begin{split}
C_{n,k}&=\frac{1}{k!}\left. \left<\frac{\left(\frac{e^{2t}+1}{2}\right)^r}{e^{\frac{1}{4}(2t)^2}}(2t)^k\right| x^n\right> \\
&=\frac{1}{k!}\left. \left<\frac{\left(\frac{e^{t}+1}{2}\right)^r}{e^{\frac{1}{4}t^2}}t^k\right| (2x)^n\right> \\
&=\frac{1}{k!}2^n \left. \left<\left(\frac{e^t+1}{2}\right)^re^{-\frac{1}{4}t^2}\right|t^k x^n \right> \\
&=2^n\binom{n}{k}\left. \left<\left(\frac{e^t+1}{2}\right)^r\right|\sum_{l=0} ^{\infty} \frac{(-1)^l}{l!4^l}t^{2l}x^{n-k}\right> \\
&=2^{n-r}\binom{n}{k}\sum_{l=0} ^{\left[\frac{n-k}{2}\right]} \frac{(-1)^l}{l!2^{2l}}(n-k)_{2l} \left. \left<\left(e^t+1\right)^r \right|x^{n-k-2l}\right>\\
&=\frac{1}{2^r} \sum_{j=0} ^r \sum_{l=0} ^{\left[\frac{n-k}{2}\right]} \frac{\binom{n}{k}\binom{r}{j}2^k(-1)^l(n-k)!}{l!(n-k-2l)!}(2j)^{n-k-2l}.
\end{split}
\end{equation}
Therefore, \eqref{23} and \eqref{24}, we obtain the following theorem.
\begin{theorem}
For $n \geq 0$, we have
\begin{equation*}
H_n(x)=\frac{1}{2^r}\sum_{k=0} ^n \left\{\sum_{j=0} ^r \sum_{l=0} ^{\left[\frac{n-k}{2}\right]} \frac{\binom{n}{k}\binom{r}{j}2^k(-1)^l(n-k)!(2j)^{n-k-2l}}{l!(n-k-2l)!}\right\}E_k ^{(r)} (x).
\end{equation*}
\end{theorem}
Note that $H_n(x)\sim\left(e^{\frac{1}{4}t^2},\frac{t}{2}\right)$. Thus, we have
\begin{equation}\label{25}
e^{\frac{1}{4}t^2}H_n(x)\sim\left(1,\frac{t}{2}\right),{\text{ and }}(2x)^n\sim\left(1,\frac{t}{2}\right).
\end{equation}
From \eqref{25}, we have
\begin{equation}\label{26}
e^{\frac{1}{4}t^2}H_n(x)=(2x)^n \Leftrightarrow H_n(x)=e^{-\frac{1}{4}t^2}(2x)^n.
\end{equation}
By \eqref{24} and \eqref{26}, we also see that
\begin{equation}\label{27}
\begin{split}
C_{n,k}&=\frac{1}{k!}\left. \left<\frac{\left(\frac{e^{2t}+1}{2}\right)^r}{e^{\frac{1}{4}(2t)^2}}(2t)^k\right| x^n\right> \\
&=\frac{1}{k!} \left. \left<\left(\frac{e^t+1}{2}\right)^rt^k\right|e^{-\frac{1}{4}t^2}(2x)^n \right> \\
&=\frac{1}{k!2^r}\left< \left. (e^t+1)^r\right|t^kH_n(x)\right>=\frac{1}{2^r}\binom{n}{k}2^k \sum_{j=0} ^r \binom{r}{j}H_{n-k} (j).
\end{split}
\end{equation}
Therefore, by \eqref{23} and \eqref{27}, we obtain the following theorem.
\begin{theorem}
For $n \geq 0$, we have
\begin{equation*}
H_n(x)=\frac{1}{2^r}\sum_{k=0} ^n\binom{n}{k}2^k\left[\sum_{j=0} ^r\binom{r}{j}H_{n-k} (j) \right]E_k ^{(r)}(x).
\end{equation*}
\end{theorem}
Let us assume that
\begin{equation}\label{28}
H_n (x)=\sum_{k=0} ^n C_{n,k}B_k ^{(r)}(x).
\end{equation}
From \eqref{13}, \eqref{17} and \eqref{28}, we have
\begin{equation}\label{29}
\begin{split}
C_{n,k}&=\frac{1}{k!}\left. \left<\frac{\left(\frac{e^{2t}-1}{2t}\right)^r}{e^{\frac{1}{4}(2t)^2}}(2t)^k\right| x^n\right> \\
&=\frac{1}{k!}\left. \left<\frac{\left(\frac{e^t-1}{t}\right)^r}{e^{\frac{1}{4}t^2}}t^k\right| (2x)^n\right> \\
&=\frac{1}{k!}\left. \left<\left(\frac{e^t-1}{t}\right)^rt^k\right|{e^{-\frac{1}{4}t^2}}(2x)^n \right>.
\end{split}
\end{equation}
From \eqref{26} and \eqref{29}, we have
\begin{equation}\label{30}
C_{n,k}=\frac{1}{k!}\left. \left<\left(\frac{e^t-1}{t}\right)^rt^k\right|H_n(x) \right>.
\end{equation}
For $r>n$, by \eqref{5} and \eqref{30}, we get
\begin{equation}
\begin{split}\label{31}
C_{n,k}&=\frac{1}{k!}\left< \left.(e^t-1)^k \right. \left| \left(\frac{e^t-1}{t} \right)^{r-k}H_n(x)\right. \right> \\
&=\frac{1}{k!}\left< \left.(e^t-1)^k \right. \left| \sum_{l=0} ^n \frac{(r-k)!}{(l+r-k)!}S_2(l+r-k,r-k)t^lH_n(x)\right. \right>\\
&=\frac{1}{k!}\sum_{l=0} ^n \frac{(r-k)!}{(l+r-k)!}S_2(l+r-k,r-k)2^l(n)_l\left<(e^t-1)^k|H_{n-l}(x)\right> \\
&=\frac{1}{k!}\sum_{l=0} ^n \frac{(r-k)!}{(l+r-k)!}S_2(l+r-k,r-k)2^l \frac{n!}{(n-l)!}\sum_{j=0} ^k\binom{k}{j}(-1)^{k-j}H_{n-l}(j) \\
&=n!\sum_{j=0} ^k \sum_{l=0} ^n \frac{(r-k)!S_2(l+r-k,r-k)(-1)^{k-j}\binom{k}{j}2^lH_{n-l}(j)}{(l+r-k)!k!(n-l)!}.
\end{split}
\end{equation}
Therefore, by \eqref{28} and \eqref{31}, we obtain the following theorem.
\begin{theorem}
For $r>n\geq0$,we have
\begin{equation*}
H_n(x)=n!\sum_{k=0} ^n \left\{\sum_{j=0} ^k \sum_{l=0} ^n \frac{(r-k)!S_2(l+r-k,r-k)(-1)^{k-j}\binom{k}{j}2^lH_{n-l}(j)}{(l+r-k)!k!(n-l)!}\right\}B_k ^{(r)}(x).
\end{equation*}
\end{theorem}
Let us assume that $r \leq n$. For $0 \leq k<r$, by \eqref{31}, we get
\begin{equation}\label{32}
C_{n,k}=n!\sum_{j=0} ^k \sum_{l=0} ^n \frac{(r-k)!S_2(l+r-k,r-k)(-1)^{k-j}\binom{k}{j}2^lH_{n-l}(j)}{(l+r-k)!k!(n-l)!}.
\end{equation}
For $r \leq k\leq n$, by \eqref{30}, we get
\begin{equation}\label{33}
\begin{split}
C_{n,k}&=\frac{1}{k!}\sum_{j=0} ^r \binom{r}{j}(-1)^{r-j}\left<e^{jt}|D^{k-r}H_n(x)\right>\\
&=\frac{2^{k-r}n!}{k!(n-k+r)!}\sum_{j=0} ^r\binom{r}{j}(-1)^{r-j}H_{n-k+r} (j).
\end{split}
\end{equation}
Therefore, by \eqref{28}, \eqref{32} and \eqref{33}, we obtain the following theorem.
\begin{theorem}
For $n \geq r$, we have
\begin{equation*}
\begin{split}
H_n(x)=&n!\sum_{k=0} ^{r-1}\left\{\sum_{j=0} ^k\sum_{l=0} ^n \frac{(r-k)!S_2(l+r-k,r-k)(-1)^{k-j}\binom{k}{j}2^lH_{n-l}(j)}{(l+r-k)!k!(n-l)!}\right\}B_k ^{(r)}(x)\\
&+n!\sum_{k=r} ^n \left\{\sum_{j=0} ^r \frac{(-1)^{r-j}\binom{r}{j}2^{k-r}H_{n-k+r}(j)}{k!(n-k+r)!}\right\}B_k ^{(r)}(x).
\end{split}
\end{equation*}
\end{theorem}
Let us assume that
\begin{equation}\label{33}
H_n (x)=\sum_{k=0} ^n C_{n,k}H_k ^{(r)}(x|\lambda).
\end{equation}
Then, by \eqref{13}, \eqref{20} and \eqref{33}, we get
\begin{equation}\label{34}
\begin{split}
C_{n,k}&=\frac{1}{k!}\left. \left<\frac{\left(\frac{e^{2t}-\lambda}{1-\lambda}\right)^r}{e^{\frac{1}{4}(2t)^2}}(2t)^k\right| x^n\right> \\
&=\frac{1}{k!}\left. \left<\frac{\left(\frac{e^t-\lambda}{1-\lambda}\right)^r}{e^{\frac{1}{4}t^2}}t^k\right| (2x)^n\right> \\
&=\frac{1}{k!} \left. \left<\left(\frac{e^t-\lambda}{1-\lambda}\right)^rt^k\right| e^{-\frac{1}{4}t^2}(2x)^n\right>.
\end{split}
\end{equation}
By \eqref{26} and \eqref{34}, we get
\begin{equation}\label{35}
\begin{split}
C_{n,k}&=\frac{1}{k!}\left. \left<\left(\frac{e^t-\lambda}{1-\lambda}\right)^rt^k\right| H_n(x)\right>=\frac{1}{k!(1-\lambda)^r}\left< \left. (e^t-\lambda)^r\right|t^kH_n(x)\right> \\
&=\frac{\binom{n}{k}2^k}{(1-\lambda)^r} \sum_{j=0} ^r\binom{r}{j}(-\lambda)^{r-j}\left<e^{jt}|H_{n-k} (x)\right> \\
&=\frac{\binom{n}{k}2^k}{(1-\lambda)^r}\sum_{j=0} ^r\binom{r}{j}(-\lambda)^{r-j}H_{n-k}(j).
\end{split}
\end{equation}
Therefore, by \eqref{33} and \eqref{35}, we obtain the following theorem.
\begin{theorem}
For $n \geq 0$, we have
\begin{equation*}
H_n (x)=\frac{1}{(1-\lambda)^r}\sum_{k=0} ^n \binom{n}{k}2^k\left[\sum_{j=0} ^r\binom{r}{j}(-\lambda)^{r-j}H_{n-k}(j)\right]H_k ^{(r)}(x|\lambda).
\end{equation*}
\end{theorem}
{\scshape Remark.} From \eqref{34}, we have
\begin{equation}\label{36}
\begin{split}
C_{n,k}&=\frac{1}{k!}\left. \left<\frac{\left(\frac{e^t-\lambda}{1-\lambda}\right)^r}{e^{\frac{1}{4}t^2}}t^k\right| (2x)^n\right>=\frac{2^n}{k!} \left. \left<\left(\frac{e^t-\lambda}{1-\lambda}\right)^re^{-\frac{1}{4}t^2}\right| t^k x^n\right>\\
&=\frac{(n)_k}{k!}2^n \sum_{l=0} ^{\left[\frac{n-k}{2}\right]}\frac{(-1)^l}{l!4^l}\left. \left<\left(\frac{e^t-\lambda}{1-\lambda}\right)^r\right| t^{2l}x^{n-k}\right> \\
&=\frac{\binom{n}{k}2^n}{(1-\lambda)^r}\sum_{l=0} ^{\left[\frac{n-k}{2}\right]}\frac{(-1)^l}{l!2^{2l}}(n-k)_{2l}\left. \left<\left(e^t-\lambda\right)^r\right| x^{n-k-2l}\right> \\
&=\frac{1}{(1-\lambda)^r}\sum_{j=0} ^r \sum_{l=0} ^{\left[\frac{n-k}{2}\right]}
\frac{\binom{n}{k}\binom{r}{j}2^k(-1)^l(-\lambda)^{r-j}(n-k)!(2j)^{n-k-2l}}{l!(n-k-2l)!}.
\end{split}
\end{equation}
Thus, by \eqref{33} and \eqref{36}, we get
\begin{equation*}
H_n(x)=\frac{1}{(1-\lambda)^r}\sum_{k=0} ^n \left\{\sum_{j=0} ^r \sum_{l=0} ^{\left[\frac{n-k}{2}\right]}
\frac{\binom{n}{k}\binom{r}{j}2^k(-1)^l(-\lambda)^{r-j}(n-k)!(2j)^{n-k-2l}}{l!(n-k-2l)!}\right\}H_k ^{(r)}(x|\lambda).
\end{equation*}

\end{document}